\documentclass[12pt]{article}
\usepackage{amsfonts,mathrsfs}
\DeclareMathSymbol{\gtrless}{\mathrel}{AMSa}{"3F}

\def\hf{{\textstyle{1 \over 2}}}
\def\defi{\stackrel{\rm def}{=}}
\def\si{\!\!\! &}
\def\sf{& \!\!\!}

\title{Zeta functions over zeros\\
of general zeta and $L$-functions}
\author{{\bf Andr\'e Voros}\footnote{Also at:
Institut de Math\'ematiques de Jussieu--Chevaleret (CNRS UMR 7586),
Universit\'e Paris 7, F-75251 Paris Cedex 05 (France).}\\
\\
CEA, Service de Physique Th\'eorique de Saclay\\
CNRS URA 2306\\
F-91191 Gif-sur-Yvette CEDEX (France)\\
E-mail : {\tt voros@spht.saclay.cea.fr}}

\begin{document}
\maketitle

We describe in detail three distinct families of generalized zeta functions
built over the (nontrivial) zeros 
of a rather general arithmetic zeta or $L$-function,
extending the scope of two earlier works that treated the Riemann zeros only.
Explicit properties are also displayed more clearly than before.
Several tables of formulae cover the simplest concrete cases:
$L$-functions for real primitive Dirichlet characters,
and Dedekind zeta functions.

\section{Generalities}

This text is a partial expansion of our oral presentation
(which surveyed an earlier paper \cite{Vz} on zeta functions
built over the Riemann zeros). Here we fully develop the argument of Sec.~5.5
therein, where we indicated how to change the Riemann zeta function $\zeta(x)$
into a more general ``primary function" (i.e., the function providing the zeros
on which the newer zeta functions are built, in the language of \cite{CHA}).
We also incorporate and extend subsequent work \cite{VZ}.
Accordingly, we can now reword the formalism to accommodate 
{\sl three\/} distinct kinds of generalized zeta functions
built over the (nontrivial) zeros 
of a fairly arbitrary (number-theoretic) zeta or $L$-function.
The resulting explicit special values are presented in (seven) Tables.
\medskip

Earlier {\sl explicit\/} descriptions of such zeta functions, 
over zeros generalizing those of $\zeta(s)$, hardly exist in the literature. 
We set apart the zeros of Selberg zeta functions: these zeros correspond to 
eigenvalues of Laplacians, and zeta functions over them
have been analyzed by spectral methods \cite{R,S,CV,Vs}
(in the cocompact case, they are particular instances of 
Minakshisundaram--Pleijel zeta functions). Otherwise,
only Dedekind (with Selberg) zeta functions get some mention \cite{K,HKW};
already the works on $L$-series \cite{Ka,JL,DS,I} 
do not discuss zeta functions {\sl per se\/} over the zeros,
but exclusively Cram\'er functions 
$V(t) \approx \sum_{\{ \rm Im \, \rho \ge 0 \}} {\rm e}^{\rho t}$,
which are somewhat related but this connection is not covered there either).
More general references are listed in greater detail in our previous articles 
\cite{Vz,VZ}.
\medskip

As for notations, we basically follow \cite{AS,B,Da}:
\begin{equation}
\label{NOT}
\matrix{ B_n \! \si : \sf \mbox{Bernoulli numbers}; \qquad \hfill 
B_n(\cdot) \! \si : \sf \mbox{Bernoulli polynomials}; \cr
E_n \! \si : \sf \mbox{Euler numbers}; \hfill 
\gamma \! \si : \sf \mbox{Euler's constant}; \hfill \cr}
\end{equation}
\vskip -.5cm
\begin{eqnarray}
\label{BET}
\beta(s) \si \defi \sf  \sum\limits_{k=0}^\infty (-1)^k (2k+1)^{-s}
\mbox{ : a specific Dirichlet $L$-function, with} \\
\!\!\! \beta(-n) = \hf E_n \!\!\!\!\!\!  && \!\!\!\!\!\!  (n \in \mathbb N) 
\ \ (e.g.,\ \beta(0)=\hf), \ \ \beta'(0) = 
-\textstyle{3 \over 2} \log 2 -\log \pi + 2 \log \Gamma ({1 \over 4}) ; 
\nonumber \\
\label{HZ}
\zeta(s,w) \si \defi \sf \sum\limits_{k=0}^\infty (k+w)^{-s}
\mbox{ : the Hurwitz zeta function.}
\end{eqnarray}
For fixed $w$, $\zeta(s,w)$ has a single pole at $s=1$, simple and of residue 1,
and the special values (\cite{B}, Sect.~1.10)
\begin{eqnarray}
\label{HN}
\zeta(-n,w) \si = \sf -B_{n+1}(w)/(n \!+\! 1) \ \ (n \in \mathbb N),
\quad (e.g.,\ \zeta(0,w)= \hf \!-\! w) \\
\label{HP}
\!\!\!\!\!\!\!\! {\rm FP}_{s=1}\, \zeta(s,w) \si = \sf -\Gamma'(w)/\Gamma(w) 
\qquad \mbox{(FP $\defi$ finite part at a pole)} \\
\label{HD}
\zeta'(0,w) \si = \sf \log \, \bigl[ \Gamma(w) / (2\pi)^{1/2} \bigr] 
\end{eqnarray}
(upon parametric zeta functions as in (\ref{HD}) and (\ref{Z1Def}--\ref{Z3Def}) 
below, $'$ will always mean differentiation with respect to 
the {\sl first\/} variable: the exponent, $s$ or $\sigma$).

Our notations (otherwise consistent with \cite{VZ}) are now generic:
objects relative to the primary zeta or $L$-function (to be denoted $L(x)$) 
will obviously depend on it, usually without mention. 

\subsection{The primary functions $L(x)$}

For the sake of definiteness, we confine interest here to situations 
still fairly close to the ``Riemann case" $L(x) = \zeta(x)$;
namely, to primary functions $L(x)$:
\begin{eqnarray}
\label{MR}
\mbox{meromorphic in ${\mathbb C}$
with at most one pole: $x=1$ (of order $q=0$ or 1);} && \quad 
\end{eqnarray}
nonvanishing in $\{ {\rm Re\,} x >1 \}$, with the normalized asymptotic behavior
\begin{equation}
\label{LAs}
(\log L)^{(n)}(x) = O(x^{-\infty}) \quad \mbox{for } x \to +\infty \quad 
(\forall n \in {\mathbb N});
\end{equation}
obeying a functional equation of the same type as $\zeta(x)$,
\begin{eqnarray}
\label{FE}
\Xi(x) \si \equiv \sf \Xi (1-x),
\qquad \Xi(x) \defi {\mathbf G}^{-1}(x) (x-1)^q L(x) ,
\end{eqnarray}
where both $\Xi(x)$ and ${\mathbf G}(x)$ are {\sl entire functions of order\/} 
$\mu_0=1$.
${\mathbf G}(x)$ is a {\sl fully explicit\/} factor 
(containing inverse Gamma factors) supplying $L(x)$ 
with explicit (``trivial") zeros $x_k \in\{ {\rm Re\,} x \le 0 \}$.
$\Xi(x)$ supplies the remaining (nontrivial) zeros of $L(x)$, 
which lie in symmetrical pairs within the strip $\{ 0< {\rm Re\,} x <1 \}$
and can be labeled as
\begin{equation}
\label{ZER}
\{\rho = \hf \pm {\rm i} \tau_k \}_{k=1,2,\ldots}, \qquad
\mbox{with ${\rm Re\,} \tau_k > 0$ and non-decreasing} 
\end{equation}
(for simplicity, we exclude the exceptional occurrence of real zeros 
$\rho$ here). 
Note: all zeros are systematically counted with multiplicities.
\smallskip

A special notation will be sometimes useful for this Taylor series at $x=1$,
\begin{equation}
\label{CUM}
\log \, \bigl[ (x-1)^q L(x) \bigr] \equiv 
\sum_{n=0}^\infty {(-1)^{n-1} \over n!} \, g_n^{\rm c} \{ L \} \, (x-1)^n .
\end{equation}
The coefficients $ g_n^{\rm c} \{ L \}$ qualify as 
{\sl generalized Stieltjes cumulants\/}:
in the case $L(x)=\zeta(x)$, $q=1$, then $g_n^{\rm c}=n \gamma_{n-1}^{\rm c}$, 
where $\{ \gamma_{n-1}^{\rm c} \}$ constitutes a cumulant sequence 
for the {\sl Stieltjes constants\/} $\gamma_{n-1}$ 
(in our notation \cite{Vz}; cf. also the $\eta_{n-1}$ in \cite{BL}, Sec.~4;
here we switch to a normalization we find more natural).

\subsection{Three zeta families}

We can then describe three (inequivalent) parametric zeta functions over 
the nontrivial zeros $\{ \rho \}$ of a generic $L$ 
satisfying the above conditions (\ref{MR}--\ref{FE}):
\begin{eqnarray}
\label{Z1Def}
{\mathscr Z}(s,x) \si \defi \sf \sum_\rho (x - \rho)^{-s} \equiv 
\sum_\rho (\rho + x - 1)^{-s} \qquad ({\rm Re\,} s >1) \\
(\mbox{definable for} \si\sf \!\!\!\!\!\!
(x - \rho) \notin {\mathbb R}^- \ (\forall \rho),
\mbox{ cf. \cite{Dn,SS,VZ} for the Riemann case}); \nonumber \\
\label{Z2Def}
{\mathcal Z}(\sigma,v) \si \defi \sf \sum_{k=1}^\infty ({\tau_k}^2+v)^{-\sigma}
\qquad \qquad \qquad \qquad \quad \ \ ({\rm Re\,} \sigma > \hf) \\
(\mbox{definable for} \si\sf \!\!\!\!\!\!
({\tau_k}^2+v) \notin {\mathbb R}^- \ (\forall k),
\mbox{ cf. \cite{G1,Dl,K,Vz} for the Riemann case}); \nonumber \\
\label{Z3Def}
{\mathfrak Z}(\sigma,y) \si \defi \sf \sum_{k=1}^\infty ({\tau_k}+y)^{-2\sigma}
\qquad \qquad \qquad \qquad \quad \ \ ({\rm Re\,} \sigma > \hf) \\
(\mbox{definable for} \si\sf \!\!\!\!\!\!
({\tau_k}+y) \notin {\mathbb R}^- \ (\forall k), 
\mbox{ cf. \cite{Vz,HKW} for the Riemann case}) , \nonumber 
\end{eqnarray}
with shorthand names for interesting special parameter values:
\begin{equation}
{\mathscr Z}(s) \!\defi\! {\mathscr Z}(s,1) , \quad 
{\mathcal Z}(\sigma) \!\defi\! {\mathcal Z}(\sigma,0) \quad
\mbox{(and, in \cite{Vz}: } 
Z(\sigma) \!\defi\! {\mathcal Z}(\sigma,{\textstyle{1 \over 4}})) .
\end{equation}

Each family is a generalized zeta function {\sl \`a la\/} Hurwitz 
(cf. eq.(\ref{HZ}));
its analytic structure mainly interests us in the exponent variable
($s$ or $\sigma$); the translation variable ($x$, $v$, or $y$) serves
to generate a parametric family in its specified range of values.

The families $\{ {\mathscr Z} \}$ and $\{ {\mathcal Z} \}$ 
share a single function, through the relation
\begin{equation}
\label{OV}
{\mathcal Z}(\sigma)\ (\equiv {\mathcal Z}(\sigma,0)) \equiv 
(2\cos\pi\sigma)^{-1} {\mathscr Z}(2\sigma,\hf).
\end{equation}
The family $\{ {\mathscr Z} \}$ can be generated from the family 
$\{ {\mathfrak Z} \}$ (but not vice-versa), thanks to the identity
(for, e.g., ${\rm Re\,} t >0$):
\begin{equation}
\label{Z2X}
{\mathscr Z} \bigl( s,\hf+t \bigr) \equiv 
\bigl[ {\rm e}^{{\rm i}\pi s/2} {\mathfrak Z} \bigl( \hf s, {\rm i} t \bigr) +
 {\rm e}^{-{\rm i}\pi s/2} {\mathfrak Z} \bigl( \hf s, -{\rm i} t \bigr) \bigr] .
\end{equation}

The families $\{ {\mathscr Z} \}$ and $\{ {\mathcal Z} \}$ are built by 
summations over all zeros $\bigl( \hf \pm {\rm i} \tau_k \bigr) $ symmetrically;
due to resulting cancellations, they will be better behaved overall 
than the third family $\{ {\mathfrak Z} \}$ 
based on the zeros with only one sign 
(this type is dubbed ``half zeta function" in \cite{HKW}).
Indeed, the first two families are formally expressible by
``explicit formulae" {\sl \`a la\/} Weil with suitably chosen test functions 
\cite{H}; however, these formulae strictly {\sl diverge\/} 
outside of clear-cut parameter domains: 
$\{ {\rm Re\,} x >1 \}$ for $\{ {\mathscr Z}(s,x) \}$,
resp. $\{ {\rm Re\,} v^{1/2} > \hf \}$ for $\{ {\mathcal Z}(\sigma,v) \}$, 
and that excludes the most interesting special cases for us: 
$x=\hf$ and 1, resp. $v=0$ and ${1 \over 4}$.
So, better adapted analytical schemes are still needed.
On the basis of all the algorithms used in \cite{Vz,VZ}, 
we see as the most efficient approach to fully describe 
the family $\{ {\mathscr Z} \}$ first, then to derive the remaining 
algebraic properties through expansion formulae in the auxiliary parameter,
and the transcendental properties from zeta-regularized factorizations of $L$.
Thus, within the same basic framework as before, 
we will get broader results in fewer steps.
Only for the full mathematical justifications must we still invoke 
\cite{Vz,VZ}.

\subsection{Range of application, and examples}

Some of the restrictions made above are just convenient to keep the paper short 
and close to concrete cases, and can probably be weakened.
For instance, zeta functions over zeros of Selberg zeta functions
have yielded results comparable to the Riemann case earlier
\cite{R,S,K,CV}, while they correspond to $\mu_0=2$ 
(${\mathbf G}$ contains a Barnes $G$-function), $q=-1$. 
Other extensions (e.g., to Hecke $L$-functions, 
as achieved upon their Cram\'er functions \cite{I}) are equally conceivable.
Currently, the assumptions are meant to closely fit two basic classes 
of examples (the most explicit properties of their zeta functions
over their zeros will be listed in a concluding section).
\medskip

{\bf - Dirichlet $L$-functions for real primitive Dirichlet characters:}
a {\sl Dirichlet $L$-function\/} is associated with a character $\chi$ of
a multiplicative group of integers mod $d$ ($d \in {\mathbb N}^\ast$ 
is called the {\sl modulus\/} or {\sl conductor\/}), as \cite{B,Da}
\begin{equation}
\label{KH}
L_\chi(x) \defi \sum_{k=1}^\infty \chi(k) \, k^{-x} \equiv 
\! \prod_{\rm \{ primes \} } \! \bigl( 1- \chi(p) \, p^{-x} \bigr) ^{-1} 
\qquad ({\rm Re\,} x >1) .
\end{equation}
Such a character is either even or odd, with a parity index $a$ defined by
\begin{equation}
\label{A01}
a = 0 \mbox{ or 1, \quad according to } \chi(-1)=(-1)^a . 
\end{equation}
$L_\chi(x)$ always satisfies conditions (\ref{MR}--\ref{LAs}) above.
We now restrict to {\sl primitive\/} characters \cite{Da},
and $d>1$ (to exclude the case $\chi \equiv 1$, $L_\chi(x) \equiv \zeta(x)$,
which more readily fits our next class of examples); then,
$L_\chi(x)$ is {\sl entire\/}, and the following functional equation holds:
\begin{eqnarray}
\label{FC}
\Xi_\chi(x) \si \equiv \sf W_\chi \, \Xi_{\overline \chi}(1-x) , \\
\label{XL}
{\rm with} \quad \Xi_\chi(x) \si \defi \sf 
(d/ \pi)^{x/2} \Gamma \bigl( \hf(x+a) \bigr) \, L_\chi(x) , \\
W_\chi \si = \sf (-{\rm i})^a d^{-1/2} 
\! \sum_{n \bmod d} \! \chi(n) {\rm e}^{2\pi{\rm i} n/d}
\end{eqnarray}
(the latter sum is called the {\sl Gaussian sum\/} for $\chi$).
The {\sl real\/} (${\overline \chi} = \chi$) primitive characters (mod $d$) 
are given by Kronecker symbols for quadratic number fields of discriminant 
$\pm d$; their Gaussian sums are explicitly known (\cite{He} thm~164),
implying $W_\chi \equiv +1$ always;
by way of consequence, the functional equation for their $L$-functions 
reduces to eq.(\ref{FE}), with 
\begin{equation}
\label{LV}
q \equiv 0; \qquad 
{\mathbf G}(x) \equiv (\pi /d)^{x/2} / \Gamma \bigl( \hf(x+a) \bigr) ,
\quad a= \biggl\{ \matrix{ 0 & \mbox{ for $\chi$ even} \cr 
1 & \mbox{ for $\chi$ odd.} \hfill \cr} \Biggr.
\end{equation}

{\bf - Dedekind $\zeta$-functions:} for any algebraic number field $K$,
its zeta function is defined as \cite{He}
\begin{equation}
\label{DZ}
\zeta_K(x) \defi \sum_{\mathfrak a} N({\mathfrak a})^{-x} \equiv 
\prod_{\mathfrak p} \bigl( 1-N({\mathfrak p})^{-x} \bigr) ^{-1} 
\qquad ({\rm Re\,} x >1)
\end{equation}
where $\mathfrak a$ (resp. $\mathfrak p$) runs over all integral (resp. prime)
ideals of $K$ and $N({\mathfrak a})$ is the norm of $\mathfrak a$.
Then $L(x)=\zeta_K(x)$ satisfies all conditions (\ref{MR}--\ref{FE}) above, with
\begin{equation}
\label{DV}
q \equiv 1, \qquad {\mathbf G}(x) \equiv 
{ \bigl( 4^{r_2} \pi^{n_K} / |d_K| \bigr)^{x/2} \over 
x \Gamma (x/2)^{r_1} \Gamma (x)^{r_2} } \, ,
\end{equation}
where $r_1$ (resp. $2r_2$) is the number 
of real (resp. complex) conjugate fields of $K$, 
$n_K \ (\equiv r_1+2r_2)$ the degree of $K$, 
and $d_K\ (\gtrless 0)$ its discriminant (\cite{He} Sec.~42, \cite{St}).
For $K = {\mathbb Q}$, having $r_1=1$, $r_2=0$ and $d_K=1$,
Riemann's $\zeta (s)$ and its classic functional equation are recovered,
with $q=1$ and 
${\mathbf G}(x) \equiv \pi^{x/2} \bigl[ x \Gamma(x/2) \bigr] ^{-1}$ 
in eq.(\ref{FE}).

\section{The first family $\{ {\mathscr Z} (s,x) \}$}

\subsection{The zeta function over the trivial zeros}

We know \cite{Vz,VZ} that a key role must be played by the zeta function 
wholly analogous to ${\mathscr Z}(s,x)$ but built on the trivial zeros of $L(x)$,
\begin{equation}
\label{SDef}
{\mathbf Z}(s,x) \defi \sum_k (x-x_k)^{-s} \qquad ({\rm Re\,} s >1) 
\end{equation}
(which we call the shadow zeta function of ${\mathscr Z}(s,x)$).
Here this function (and all quantities referring to it) is to be viewed
as fully explicitly known, just as the trivial zeros themselves are
(in our concrete examples, ${\mathbf Z}(s,x)$ will be expressed 
in terms of the Hurwitz zeta function (\ref{HZ})). 
It is also necessary to relate ${\mathbf Z}(s,x)$ and ${\mathbf G}(x)$ 
as fully as possible.
We specially want the results {\sl to include formulae for the special values\/}
${\rm FP}_{s=1} {\mathbf Z} (s,x)$ {\sl and\/} 
${\rm FP}_{s=1} {\mathscr Z} (s,x)$,
{\sl resp.\/} ${\mathbf Z}' (0,x)$ {\sl and\/} ${\mathscr Z}' (0,x)$,
because comparable formulae for the Hurwitz zeta function are quite important:
eqs.(\ref{HP}), resp. (\ref{HD}). 
Then the original normalization of the trivial factor is not fully adequate:
it is better to rewrite $L(x)$ as a product of either Weierstrass-like factors
(as in \cite{Vz}) or zeta-regularized factors (as in \cite{VZ});
we follow the latter course here. 
\medskip

We recall some off-the-shelf rules on zeta-regularization for
suitable infinite products of order $\mu_0=1$, of the form
$\Delta (x) = {\rm e}^{B_1 x + B_0} \prod_k [(1-x/y_k) {\rm e}^{x/y_k}] $
(\cite{VZ} and refs. therein;
the rules will be valid for ${\mathbf G}(x)$ and $\Xi(x)$).
Actually, eqs.(\ref{STD}) and (\ref{ZLM}) will also serve for $\mu_0<1$ 
in Sec.~3, and (only they) are worded for general $0< \mu_0 \le 1$.
A key further requirement is an {\sl asymptotic\/} (``generalized Stirling" \cite{JL}) 
{\sl expansion for\/} $\log \Delta(x)$ $(x \to +\infty)$, of the form
\begin{equation}
\label{STD}
\log \Delta(x) \sim \tilde a_1 x (\log x -1) + b_1 x + \tilde a_0 \log x + b_0 
\ \ \Bigl[ + \!\!\! \sum_{ \{ \mu_k \} \setminus \{ 0,1\} } a_{\mu_k} x^{\mu_k}
\Bigr]
\end{equation}
for some sequence $(1 \ge) \ \mu_0 > \mu_1 > \cdots > \mu_n \downarrow -\infty$,
and indefinitely differentiable term by term. 
Here, the terms designated by coefficients $a_\mu$ or $\tilde a_\mu$ are 
those {\sl \b allowed\/} in a {\sl zeta-regularized\/} product;
any {\sl extra\/} terms with a pure $x^1$ or $x^0$ dependence 
are {\sl \b banned\/} (they read as $b_n x^n$ here). 
A generalized zeta function is also introduced, as 
$ Z(s,x) = \sum_k (x-y_k)^{-s} $ if ${\rm Re\,} s > \mu_0$. Then:
\smallskip

\noindent - the zeta-regularized form $D(x)$ of a product $\Delta (x)$ 
is explicitly obtained just by removing any ``banned" portion
present in the large-$x$ expansion of $\log \Delta(x)$: specifically, here,
\begin{equation}
\label{DZR}
D (x) \defi {\rm e}^{-Z'(0,x)} \equiv {\rm e}^{-(b_1 x + b_0)} \Delta (x) ;
\end{equation}
- the logarithmic derivatives of the zeta-regularized product yield
\begin{eqnarray}
\label{ZLD}
(\log D)' (x) \si \equiv \sf {\rm FP}_{s=1} Z (s,x), \\
\label{ZLM}
{(-1)^{m-1} \over (m \!-\! 1)!} (\log D)^{(m)} (x) \si \equiv \sf Z(m,x) 
\qquad \mbox{for integer } m > \mu_0 
\end{eqnarray}
(we will also use this result once with $\mu_0=\hf$, in Sec.~3);

\noindent - the results carry over to non-integer $s$ as Mellin-transform formulae:
\begin{equation}
\label{MEL}
Z(s,x) = {\sin \pi s \over \pi (1 \!-\! s)} I(s,x), \quad  
I(s,x) \defi \! \int_0^\infty \! 
Z(2,x+y) \, y^{1-s} \,{\rm d} y \quad (1 < {\rm Re\,} s < 2) ;
\end{equation}
then, $I(s,x)$ extends to a meromorphic function in the whole $s$-plane
through repeated integrations by parts, and its polar structure is 
fully encoded in the $(y \to +\infty)$ expansion of $Z(2,x+y)$,
itself computable (see example next).
\medskip

We now specialize the above results first to the trivial factor 
${\mathbf G}(x)$.
Consisting mainly of inverse Gamma factors, ${\mathbf G}(x)$ has 
a {\sl computable\/} Stirling expansion (for $x \to +\infty$) 
which can be reorganized in the form 
\begin{equation}
\label{ST}
-\log {\mathbf G}(x) \sim \tilde a_1 x (\log x -1) + b_1 x
+ \tilde a_0 \log x + b_0 \ + \sum_{n=1}^\infty a_{-n} x^{-n}
\end{equation}
(and which also governs $[\log \Xi (x) -q \log (x-1)]$, 
by eqs.(\ref{LAs}) and (\ref{FE})). Eq.(\ref{DZR}) then implies
that the zeta-regularized forms for ${\mathbf G}(x)$ and $\Xi(x)$ are
\begin{eqnarray}
\label{SZR}
{\mathbf D} (x) \si \defi \sf {\rm e}^{-{\mathbf Z}'(0,x)} \equiv 
{\rm e}^{+b_1 x + b_0} {\mathbf G}(x) , \\
\label{ZR}
{\mathscr D} (x) \si \defi \sf {\rm e}^{-{\mathscr Z}'(0,x)} \equiv 
{\rm e}^{-b_1 x - b_0} \Xi(x) ,
\end{eqnarray}
which in turn entail a zeta-regularized decomposition of $L(x)$, as
\begin{equation}
\label{zf1}
(x-1)^q L(x) \equiv {\mathbf D} (x) {\mathscr D} (x) .
\end{equation}
Concretely here (using eq.(\ref{SZR}), and $\mu_0=1$),
eqs.(\ref{ZLD}) and (\ref{ZLM}) translate to
\begin{eqnarray}
\label{ZG}
{\rm FP}_{s=1} {\mathbf Z} (s,x) \si \equiv \sf (\log {\mathbf G})'(x) +b_1, \\
\label{ZGM}
{\mathbf Z}(m,x) \si \equiv \sf 
{(-1)^{m-1} \over (m \!-\! 1)!} (\log {\mathbf G})^{(m)} (x) 
 \quad \mbox{for } m=2,3,\ldots .
\end{eqnarray}
The substitution of the Stirling series (\ref{ST}) into eq.(\ref{ZGM}) 
with $m=2$ leads to the $(y \to +\infty)$ expansion of ${\mathbf Z}(2,x+y)$ 
in the simple form $\sum_{n \ge -1} c_{-n}(x) y^{-n-2}$
(the $c_{-n}(x)$ are polynomials). It follows that eq.(\ref{MEL}), 
written for $Z={\mathbf Z}$, yields an $I(s,x)$ with poles at $s=-n$,
$n=-1,0,1,2,\ldots$, all simple and of residues $c_{-n}(x)$. 
As consequences for ${\mathbf Z}(s,x)$, restated in fully explicit form: 
\smallskip

\noindent - ${\mathbf Z}(s,x)$ extends to a meromorphic function 
in the whole $s$-plane, with
\begin{equation}
\label{SP1}
\mbox{the single pole $s=1$, simple, of residue $\tilde a_1$ \quad
(independent of $x$).}
\end{equation}
- the values ${\mathbf Z}(-n,x),\ n \in {\mathbb N}$ are given by 
{\sl closed polynomial formulae\/} (``trace identities" in a spectral setting),
as
\begin{equation}
\label{TI}
\matrix{ \displaystyle {\mathbf Z}(-n,x) = 
- {\tilde a_1 \over n \!+\! 1} \, x^{n+1} - \tilde a_0 x^n
+ n \sum_{j=1}^n (-1)^j {n \!-\! 1 \choose j \!-\! 1} a_{-j} x^{n-j} , \cr
e.g., \quad {\mathbf Z}(0,x) = -\tilde a_1 x - \tilde a_0 \, . \cr}
\end{equation}

Then, the same formulae overall hold with 
${\mathbf D}$, ${\mathbf Z}$, ${\mathbf G}$
replaced by the (less explicit) ${\mathscr D}$, ${\mathscr Z}$, $\Xi$ 
respectively (and with suitably changed coefficients).

\subsection{The main result}

As a generalization of eq.(42) in \cite{VZ}, ${\mathscr Z} (s,x)$ admits 
an integral representation valid in the half-plane $\{ {\rm Re\,} s <1 \}$ and 
for any eligible value of the parameter $x$ avoiding the cut $(-\infty,+1]$:
\begin{eqnarray}
\label{ZC1}
{\mathscr Z}(s,x) \si = \sf -{\mathbf Z}(s,x) + {q \over (x-1)^s}
+ {\sin \pi s \over \pi} {\mathscr J}\! (s,x),\\
\label{ZRP}
{\mathscr J}\! (s,x) \si \defi \sf 
\int_0^\infty {L' \over L} (x+y) \, y^{-s} \,{\rm d} y \qquad ({\rm Re\,} s <1);
\end{eqnarray}
here, $(x-1)^s$ is given its standard determination for 
$x \in {\mathbb C} \setminus (-\infty,+1]$; 
its discontinuity across the real axis, as well as those of $-{\mathbf Z}(s,x)$,
are precisely cancelled through corresponding jumps of ${\mathscr J}\! (s,x)$
so that only the {\sl nontrivial\/} zeros of $L(x)$ induce 
genuine $x$-plane discontinuities in ${\mathscr Z}(s,x)$
(this can be checked by comparing computations of the right-hand side with
small imaginary parts ($\pm {\rm i} 0$) added to $x$).

Eq.(\ref{ZC1}) easily takes real forms; 
a simple one valid for $x>0$ (at least) is
\begin{equation}
\label{ZR1}
{\mathscr Z}(s,x)= -{\mathbf Z}(s,x) + {\sin \pi s \over \pi} \int_0^\infty 
\Bigl[ {L' \over L} (x+y) + {q \over x \!+\! y \!-\! 1} \Bigr] y^{-s}
\,{\rm d} y \quad (0<{\rm Re\,} s <1);
\end{equation}
this form only converges in a strip of the $s$-plane, 
but unlike (\ref{ZC1}), it remains well defined as $x \to +1$:
\begin{equation}
{\mathscr Z}(s,1) = -{\mathbf Z}(s,1) + {\sin \pi s \over \pi} 
\int_0^\infty \Bigl[ {L' \over L} (1+y)+{q \over y} \Bigr] y^{-s} \,{\rm d} y .
\end{equation}

The proof is now a simple application of the general formulae above: 
the second logarithmic derivative of eq.(\ref{zf1}) (or (\ref{FE})) gives
$ [ {\mathscr Z} \!+\! {\mathbf Z} ] (2,x) = q (x \!-\! 1)^{-2} -(L'/L)'(x) $
(using eq.(\ref{ZLM}) for $D={\mathbf G} \, \Xi$ and $m=2$); 
upon this and for $x \notin (-\infty,+1]$, the Mellin formula (\ref{MEL}) yields
\begin{equation}
[ {\mathscr Z} + {\mathbf Z} ] (s,x) \equiv
{\sin \pi s \over \pi (1-s)} \! \int_0^\infty  
\Bigl[ {q \over (x+y-1)^2 } - \Bigl( {L' \over L} \Bigr) '(x+y) \Bigr]
\, y^{1-s} \,{\rm d} y .
\end{equation}
Now, moreover, {\sl this\/} integral is convergent and analytic 
in $ \{ 0 < {\rm Re\,} s < 2 \} $,
the integrand being regular, and $O(1)$ for $y \to +0$, 
and (using eq.(\ref{LAs})) $O(y^{-2})$ for $y \to +\infty$. 
Then, splitting the integral in two, 
the first term evaluates explicitly to $q(x-1)^{-s}$,
and upon restricting to $0 < {\rm Re\,} s < 1$, 
the second term transforms to $\pi^{-1} \sin \pi s {\mathscr J}\! (s,x)$ 
through an integration by parts; hence the result (\ref{ZC1}).
(The same integration by parts upon the whole integral
yields eq.(\ref{ZR1}) instead.)

\subsection{Explicit consequences for the family $\{ {\mathscr Z}(s,x) \}$}

The subsequent statements refer to $s$ as variable, with $x$ fixed. 

Eq.(\ref{ZC1}) is a true analog for ${\mathscr Z}(s,x)$ 
of the (Joncqui\`ere--Lerch) functional relation for $\zeta(s,w)$ 
(\cite{B}, Sect.~1.11 eq.(16)).
First of all, it gives an explicit one-step analytical continuation 
of ${\mathscr Z}(s,x)$ to the half-plane $\{{\rm Re\,} s<1 \}$. It also implies 
its {\sl meromorphic continuation in\/} $s$ to all of ${\mathbb C}$, 
since a Mellin transform like ${\mathscr J}\! (s,x)$ has 
a well understood meromorphic structure: repeated integrations by parts
on eq.(\ref{ZRP}) (invoking the asymptotic formula (\ref{LAs})) 
show that ${\mathscr J}\! (s,x)$ is meromorphic in the whole $s$-plane,
and has only simple poles at $s=1,2,\ldots$, with residues
\begin{equation}
\label{REZ}
{\rm Res}_{s=n} {\mathscr J}\! (s,x)
= -(\log |L|)^{(n)}(x)/(n-1)! \qquad (n=1,2,\ldots) .
\end{equation}
It then follows from eq.(\ref{ZC1}) that ${\mathscr Z}(s,x)$ precisely inherits
the polar structure of $-{\bf Z}(s,x)$, namely (cf. eq.(\ref{SP1})):
\begin{equation}
\label{ZP1}
{\mathscr Z}(s,x) \mbox{ has the single pole $s=1$, simple, 
of residue $-\tilde a_1$.}
\end{equation}

If $L(x)$ admits an Euler product, as in the examples 
(\ref{KH}) and (\ref{DZ}), then the substitution of its logarithmic derivative
into eq.(\ref{ZRP}), followed by integration term by term, yields 
an {\sl asymptotic $(s \to -\infty)$ expansion\/} for ${\mathscr J}\! (s,x)$, 
and thereby for ${\mathscr Z}(s,x)$, just as in the Riemann case 
(\cite{VZ}, eq.(52)).
\medskip

Finally, almost all the {\sl special values\/} of ${\mathscr Z}(s,x)$ 
(at integer $s$) are explicitly readable off eq.(\ref{ZC1}) 
thanks to its vanishing factor ($\sin \pi s$). 
On general grounds (\cite{Vz}, Sec.~4),
the ${\mathscr Z}(n,x)$ come out {\sl algebraically\/} for $n \in -{\mathbb N}$
(as given in Table~1, upper part),
and (together with ${\mathscr Z}'(0,x)$) {\sl transcendentally\/} 
for $n \in {\mathbb N}^\ast$, for instance:
\begin{eqnarray}
{\mathscr Z}'(0,x) \si = \sf 
-{\mathbf Z}'(0,x) -q \log (x-1) + {\mathscr J}\! (0,x) \nonumber \\
\label{ZDX}
\si = \sf b_1 x + b_0 + \log {\mathbf G}(x) - \log \,\bigl[(x-1)^q L(x) \bigr] 
\\
{\rm FP}_{s=1} {\mathscr Z}(s,x) \si = \sf -{\rm FP}_{s=1} {\mathbf Z}(s,x) 
+ {q \over x-1} - {\rm Res}_{s=1} {\mathscr J}\! (s,x), \nonumber \\
\label{ZFP}
\si = \sf -b_1 - (\log {\mathbf G})' (x) 
+ \biggl[ {q \over x \!-\! 1} + {L' \over L} (x) \biggr] 
\end{eqnarray}
(the reduced forms use eqs.(\ref{SZR}), (\ref{ZG}), (\ref{REZ}));
and likewise for ${\mathscr Z}(n,x)$, $n \ge 2$. 
Terms in brackets stay globally continuous for $x \to 1$ 
(also in eq.(\ref{ZPP}) below); 
the limiting special values will be expressed as eq.(\ref{Z1V}) later.

However, the values ${\mathscr Z}(n,x)$ for $n \in {\mathbb N}^\ast$ 
({\sl also including\/} $n=1$) emerge more simply 
by differentiating the logarithm of a symmetrical Hadamard product formula, 
$\Xi(x) \propto \prod_\rho (1-x/ \rho)$ 
(just as in the Riemann case \cite{VZ}), as
\begin{eqnarray}
\label{ZNP}
\!\!\!\!\!\!\!\!\!\!{\mathscr Z}(n,x) \si = \sf 
{(-1)^{n-1} \over (n-1)!} (\log \Xi)^{(n)}(x)
\qquad \qquad \qquad \qquad \qquad \quad \ \ (n=1,2,\ldots) \\
\label{ZPP}
\si = \sf -{\mathbf Z}(n,x) + \biggl[ {q \over (x \!-\! 1)^n} 
\!+\! {(-1)^{n-1} \over (n \!-\! 1)!} (\log |L|)^{(n)}(x) \biggr] 
\ \, (n=2,3,\ldots) .
\end{eqnarray}
(This is by far the quickest path to transcendental values,
but it altogether misses the pair we declared to be important:
${\mathscr Z}'(0,x)$ and ${\rm FP}_{s=1} {\mathscr Z}(s,x)$, 
determined above by eqs.(\ref{ZDX}) and (\ref{ZFP}) respectively.)

For $n=1$, eq.(\ref{ZPP}) cannot hold: ${\mathscr Z}(1,x)$ is finite 
(eq.(\ref{ZNP}) {\sl defines\/} it to be $\sum_\rho (x-\rho)^{-1}$ 
with the zeros ordered pairwise, as usual),
and ${\mathbf Z}(1,x)$ diverges (at the same time, both functions 
${\mathscr Z}(s,x)$, ${\mathbf Z}(s,x)$ have a pole at $s=1$~!). 
Instead, the comparison of eq.(\ref{ZNP}) at $n=1$ with eq.(\ref{ZFP})
yields a {\sl fixed anomaly\/}, or discrepancy between two natural
specifications for a finite value at $s=1$,
\begin{equation}
\label{FPV}
{\mathscr Z}(1,x) - {\rm FP}_{s=1} {\mathscr Z}(s,x) = 
(\log \, [\Xi/{\mathscr D}])'(x) = b_1 \quad \mbox{(constant)}. 
\end{equation}

Table~1 summarizes the special values obtained for ${\mathscr Z}(s,x)$ 
at general~$x$ (thus extending to general primary functions $L(x)$ formulae
previously restricted to $L(x)=\zeta(x)$, cf. Table~2 in \cite{VZ}).

\begin{table}
\centering
\begin{tabular} {ccc}
\hline \\[-8pt]
$s$ & &
${\mathscr Z}(s,x) = \sum\limits_\rho (x-\rho)^{-s} $ \\[7pt]
\hline \\[-12pt]
$-n \le 0$ & & $ -{\mathbf Z}(-n,x) + q(x-1)^n $ \\[2pt]
$0$ & & $ \tilde a_1 x + \tilde a_0 +q$ \\[2pt]
\hline \\[-12pt]
{\sl $s$-derivative at 0} & &
$ {\mathscr Z}'(0,x) = b_1 x + b_0 - \log \Xi(x) $ \\[2pt]
{\sl finite part at +1} & &
$ \displaystyle {\rm FP}_{s=1} {\mathscr Z}(s,x) = - b_1 + (\log \Xi)'(x) $ \\
$+n \ge 1 $ & & $ {\textstyle (-1)^{n-1} \over \textstyle (n-1)!} 
(\log \Xi)^{(n)}(x) $ \\[8pt]
\hline\\[-12pt]
\end{tabular}
\caption{Special values of ${\mathscr Z}(s,x)$
(upper part: algebraic, lower part: transcendental \cite{Vz})
for a general primary zeta function $L(x)$ with a pole of order $q$ (at $x=1$).
Notations: see eqs.(\ref{FE}) for $\Xi(x)$, 
(\ref{SDef}) and (\ref{TI}) for ${\mathbf Z}(-n,x)$,
(\ref{ST}) for $\tilde a_j,\ b_j$; $n$ is integer.
}
\end{table}

\medskip

The two sets of linear identities for the values ${\mathscr Z}(n,x)$ 
in the Riemann case (eqs.~(61--62) in \cite{VZ}), 
which are purely induced by the symmetry 
$(\rho \longleftrightarrow 1-\rho)$ in eq.(\ref{Z1Def}),
naturally persist here:
\begin{eqnarray}
\label{ZSN}
\!\!{\mathscr Z}(n,x) \si = \sf (-1)^n {\mathscr Z}(n,1-x) 
\qquad \qquad \qquad \qquad \mbox{for } n=1,2,\ldots;\\
\label{ZSK}
\!\!{\mathscr Z}(k,x) \si = \sf - {1 \over 2} \!\sum_{\ell=k+1}^\infty \!\!
{\ell \!-\! 1 \choose k \!-\! 1} (2x \!-\! 1)^{\ell-k} {\mathscr Z}(\ell,x) 
\quad \mbox{for each {\sl odd\/} } k \ge 1.
\end{eqnarray}

More explicit formulae can result for exceptional $x$-values 
such as $x=\hf$ and $x=1$,
which respectively correspond (via the functional equation (\ref{FE})) 
to the symmetry center of $\Xi (x)$ and to the origin in the $x$-plane
(${\mathscr Z}(s,x)$ is clearly regular on the real $x$-axis,
under our assumption $\{ \rho \} \cap {\mathbb R} = \emptyset$).
Thus, for $x=\hf$, eq.(\ref{ZNP})  simplifies to
\begin{eqnarray}
\label{Z10}
{\mathscr Z}(n,\hf) \si \equiv \sf 0 \qquad \mbox{for all $n \ge 1$ odd} \\
\label{Z1FP}
(\mbox{and} \quad 
{\rm FP}_{s=1} {\mathscr Z}(s,\hf) \si = \sf -b_1 \quad \mbox{by eq.(\ref{FPV})})
\end{eqnarray}
(in combination with eqs.(\ref{ZPP}),(\ref{ZFP}), 
these amount to the explicit specifications
$(\log |L|)^{(2m+1)}(\hf) = (\log {\mathbf G})^{(2m+1)}(\hf) 
+ 2^{2m+1}q \,(2m)! \,$,
also directly readable off the functional equation (\ref{FE})); while
for $x=1$, that same formula brings in the Taylor series (\ref{CUM}), to yield
\begin{eqnarray}
{\mathscr Z}(1,1) \si = \sf -(\log {\mathbf G})'(1) + g_1^{\rm c} \{ L \} , 
\nonumber \\
\label{Z1V}
{\mathscr Z}(n,1) \si = \sf -{\mathbf Z} (n,1) + g_n^{\rm c} \{ L \} /(n-1)! 
\qquad (n=2,3,\ldots) \\
(\mbox{also,} \quad {\mathscr Z}'(0,1) \si = \sf 
- {\mathbf Z}'(0,1) + g_0^{\rm c} \{ L \}). \nonumber 
\end{eqnarray}
Case by case, ${\mathbf Z}(s,x)$ can also be made explicit for $x=\hf$ or 1,
just as when $L(x)=\zeta(x)$ (\cite{VZ}, Sec.~3.3);
for our selected examples (\ref{KH}), (\ref{DZ}) it will reduce to combinations 
of the two fixed Dirichlet series $\zeta(s)$ and $\beta(s)$ 
(cf. eq.(\ref{BET})).
The accordingly reduced special values of ${\mathscr Z}(s,\hf)$ and 
${\mathscr Z}(s,1)\ (= {\mathscr Z}(s))$ are tabulated in the concluding Sec.~5.

\section{The second family $\{ {\mathcal Z} (\sigma,v) \}$}

The main starting tool is the relation (\ref{OV}), which allows a 1--1 transfer 
of the previous results for ${\mathscr Z} (s,\hf)$ onto one particular member, 
${\mathcal Z} (\sigma) \equiv {\mathcal Z} (\sigma,v=0)$, of the second family.

\subsection{The basic case $v=0$}

The identity (\ref{OV}) shows that ${\mathcal Z} (\sigma)$
is meromorphic in all of $\mathbb C$ with a double pole at $\sigma=\hf$,
simple poles at $\sigma=\hf-m$ ($m=1,2,\ldots$), and polar parts:
\begin{eqnarray}
{\mathcal Z}(\hf+\varepsilon) \si = \sf 
- {\rm Res}_{s=1} {\mathscr Z}(s,\hf) {1 \over 4\pi\varepsilon^2}
- {\rm FP}_{s=1} {\mathscr Z}(s,\hf) {1 \over 2\pi\varepsilon} 
\ + O(1)_{\varepsilon \to 0} \nonumber \\
\label{Z2P}
\si = \sf {\tilde a_1 \over 4\pi\varepsilon^2} + {b_1 \over 2\pi\varepsilon} 
\ + O(1)_{\varepsilon \to 0} \ ; \\
{\mathcal Z}(\hf-m+\varepsilon) \si = \sf 
{(-1)^{m+1} \over 2\pi\varepsilon} {\mathscr Z}(1-2m,\hf) 
\ + O(1)_{\varepsilon \to 0} \qquad \mbox{for } m=1,2,\ldots , \nonumber \\
\label{Z2Q}
i.e., \quad 
{\mathcal R}_m \si \defi \sf 
{\rm Res}_{\sigma={1 \over 2}-m} {\mathcal Z}(\sigma)
= {(-1)^m \over 2\pi} \bigl[ {\mathbf Z}(1-2m,\hf) +q \, 2^{1-2m} \bigr] 
\end{eqnarray}
(eq.(\ref{Z2P}) follows from eqs.(\ref{ZP1}),(\ref{Z1FP}), 
and ${\mathcal R}_m$ from using Table~1).

That relation also yields an integral representation for ${\mathcal Z} (\sigma)$
(eq.(72) in \cite{Vz} for the Riemann case), 
just by specializing eq.(\ref{ZC1}) to $x=\hf \pm {\rm i} 0$.
Finally, it delivers all the special values of ${\mathcal Z} (\sigma)$ as 
\begin{equation}
\label{Z1E}
{\mathcal Z}(m) = \hf (-1)^m {\mathscr Z}(2m,\hf) \quad 
{\rm for\ } m \in {\mathbb Z}, \qquad {\mathcal Z}'(0) = {\mathscr Z}'(0,\hf),
\end{equation}
which become fully explicit using Table 1. The problem is then to extend 
all those results to general values of the parameter $v$.

\subsection{Algebraic results for general $v$}

Our best tool here is a straightforward expansion of ${\mathcal Z}(\sigma,v)$ 
around $v=0$ \cite{Vz}, convergent for $|v| < \min_k \{ |\tau_k|^2 \} $:
\begin{equation}
\label{SHZ0}
{\mathcal Z}(\sigma,v) = \sum_{k=0}^\infty 
({\tau_k}^2)^{-\sigma} \Bigl( 1 + {v \over {\tau_k}^2} \Bigr)^{-\sigma} 
= \sum_{\ell=0}^\infty {\Gamma(1-\sigma) \over \ell !\Gamma(1-\sigma-\ell)}
{\mathcal Z}(\sigma+\ell) \, v^\ell .
\end{equation}

This first provides the meromorphic continuation of ${\mathcal Z}(\sigma,v)$ 
at fixed $v \ne 0$ to the whole complex $\sigma$-plane, 
now with {\sl double poles\/} at all $\sigma=-m+\hf,\ m \in\mathbb N$. 
More precisely, the polar part of ${\mathcal Z}(\sigma,v)$ at each pole 
only depends on a {\sl finite\/} stretch of the series (\ref{SHZ0}),
\begin{equation}
\label{HP2}
{\mathcal Z}(-m+\hf+\varepsilon,v) = \sum_{\ell=0}^m 
{\Gamma(\hf \!+\! m \!-\! \varepsilon) \over 
\ell ! \Gamma(\hf \!+\! m \!-\! \ell \!-\! \varepsilon)}
{\mathcal Z}(-m+\ell+\hf+\varepsilon) \, v^\ell
\ + O(1)_{\varepsilon \to 0} \ ;
\end{equation}
upon importing the polar structure of ${\mathcal Z}(\sigma)$ 
from eqs.(\ref{Z2P}--\ref{Z2Q}), this yields
\begin{eqnarray}
\label{Hp}
{\mathcal Z}(-m \!+\! \hf \!+\! \varepsilon,v) =
{\tilde a_1 \over 4 \pi} 
{\Gamma(m+\hf) \over m! \Gamma(\hf)} v^m \,\varepsilon^{-2}
 + {\mathcal R}_m(v) \,\varepsilon^{-1} \ + O(1)_{\varepsilon \to 0} \ , 
\qquad \qquad \\
{\mathcal R}_m(v) = -{\Gamma(m \!+\! \hf) \over m! \Gamma(\hf)} 
\Bigl[ {\tilde a_1 \over 2 \pi} \sum_{j=1}^m {1 \over 2j \!-\! 1}
- {b_1 \over 2 \pi} \Bigr] v^m 
+ \sum_{j=1}^m{\Gamma(\hf \!+\! m) \over (m \!-\! j)!\Gamma(\hf \!+\! j)}
{\mathcal R}_j \, v^{m-j} . \nonumber
\end{eqnarray}
Here, the polar part of order 2 at every $(-m+\hf)$ clearly comes from 
the single double pole of ${\mathcal Z}(\sigma)$ (at $\sigma=\hf$); 
whereas each residue 
${\mathcal R}_m(v)$ has contributions from all the residues of 
${\mathcal Z}(\sigma)$ at $\sigma = -j+\hf \ge -m+\hf$ 
(specified by eqs.(\ref{Z2P}) for $j=0$, (\ref{Z2Q}) for $j \ge 1$).
For the leading pole $s=\hf$, eq.(\ref{Hp}) boils down to eq.(\ref{Z2P}),
i.e., {\sl this\/} full polar part is independent of $v$.

For the special values ${\mathcal Z}(-m,v)$, $m \in \mathbb N$, 
the series (\ref{SHZ0}) also terminates:
\begin{equation}
\label{TIG}
{\mathcal Z}(-m,v) \equiv 
\sum_{\ell=0}^m {m \choose \ell} {\mathcal Z}(-m+\ell)\,v^\ell 
\qquad (m \in \mathbb N) ,
\end{equation}
where the values ${\mathcal Z}(-m+\ell)$ are explicit from eq.(\ref{Z1E}) 
and Table~1; e.g.,
\begin{equation}
\label{Z0V}
{\mathcal Z}(0,v) \equiv \hf \bigl[ -{\mathbf Z}(0,\hf) + q \bigr] 
\equiv \hf(\hf \tilde a_1 + \tilde a_0 +q)
\qquad (\mbox{independent of } v) .
\end{equation}

In the end, all the polar terms of ${\mathcal Z}(\sigma,v)$, 
and the special values ${\mathcal Z}(-m,v)$ ($m \in \mathbb N$), 
(listed in Table~2, upper half), are {\sl (computable) polynomials\/} in $v$.

\subsection{Transcendental values for general $v$}

Now an optimal tool is a variant of the factorization (\ref{zf1}),
using the alternative zeta-regularized factor 
\begin{equation}
{\mathcal D}(v) \defi {\rm e}^{-{\mathcal Z}'(0,v)}
\end{equation}
instead of ${\mathscr D}(x)$. 
The main point here is the replacement of $x$ by $v=(x-\hf)^2$ 
as basic variable: then, in contrast to ${\mathscr D}(x)$ of eq.(\ref{zf1}), 
this zeta-regularization of $\Xi(x)$ now {\sl preserves the symmetry\/}
$(x \leftrightarrow 1-x)$.

Rewritten in the variable $v \to + \infty$, the Stirling formula (\ref{ST})
for $[\log \Xi (x) -q \log (x-1)]$ becomes of order $\hf <1$, yielding
\begin{equation}
\log \Xi (\sqrt v + \hf) \sim
\hf \tilde a_1 v^{1 \over 2} \log v + (b_1-\tilde a_1) v^{1 \over 2} +
\hf \bigl( \hf \tilde a_1 + \tilde a_0 + q \bigr) \log v
+ ( \hf b_1 + b_0) \ [+O(v^{-{1 \over 2}})] ;
\end{equation}
as ``banned" terms (cf. Sec.~2.1), only constants ($\propto v^0$) can now occur 
($v^\mu \log v$ are allowed if $\mu \notin {\mathbb N}$, they just induce
{\sl double\/} poles in the zeta functions \cite{Vz}); here this results in
\begin{equation}
\label{EF}
{\mathcal D}(v) \equiv {\rm e}^{-(b_0 + b_1 /2)} \Xi (\hf+v^{1/2}) \equiv 
{\rm e}^{b_1 v^{1/2}} {\mathscr D} (\hf+v^{1/2})
\end{equation}
and in the modified decomposition 
(cf. eq.(41) in \cite{VZ} for the Riemann case)
\begin{equation}
\label{zf2}
(x-1)^q L(x) \equiv {\rm e}^{-b_1 (x-1/2)} {\mathbf D}(x) {\mathcal D}(v),
\qquad v \defi (x-\hf)^2 .
\end{equation}

All transcendental special values of ${\mathcal Z}(\sigma,v)$ immediately 
follow, just as in the Riemann case \cite{Vz}:
first, ${\mathcal Z}'(0,v) \equiv -\log {\mathcal D}(v)$ 
expresses in terms of $\log \Xi(\hf \pm v^{1/2})$
(and thereby, of $\log L(\hf \pm v^{1/2})$), using eq.(\ref{EF}); 
then eq.(\ref{ZLM}), now applied to ${\mathcal D}$ and ${\mathcal Z}$ 
with $v$ as variable hence $\mu_0=\hf$, yields 
\begin{equation}
{\mathcal Z}(m,v) = {(-1)^{m-1} \over (m-1)!} \,
{{\rm d}^m \over {\rm d}v^m} \log \Xi (\hf \pm v^{1/2}), \qquad m=1,2,\ldots .
\end{equation}
By the chain rule, the right-hand side must simplify to 
a finite linear combination of derivatives $(\log \Xi)^{(\ell)}(x)$ 
at $x= \hf \pm v^{1/2}$ 
(and thereby, of values ${\mathscr Z}(\ell,x)$ by eq.(\ref{ZNP})),
but it is unclear how to carry through this nonlinear change of variables 
for general $m$ and $v$ directly. 
Instead, setting $v \equiv (x-\hf)^2$ and $\rho \equiv \hf + {\rm i}\tau$ 
throughout here, we may start from the identity
\begin{equation}
\Bigl( 1-{s \over x-\rho} \Bigr) \Bigl( 1-{s \over x-1+\rho} \Bigr) \equiv
1- {s(2x-1-s) \over \tau^2+v} \, , 
\end{equation}
expand the logarithms of both sides, and identify like powers of $s$
to get a triangular sequence (for $n=1,2,\ldots$) of linear identities,
\begin{equation}
\label{LID}
{ (x \!-\! \rho)^{-n} + (x \!-\! 1 \!+\! \rho)^{-n} \over n} \equiv 
\!\!\! \sum_{0 \le \ell \le n/2} \! (-1)^\ell {n \!-\! \ell \choose \ell} 
(2x - 1)^{n-2\ell} \,
{(\tau^2 \!+\! v)^{-n+\ell} \over n-\ell} \ .
\end{equation}
Then, on the one hand, summing this over (half) the zeros yields the identities
\begin{equation}
\label{LIZ}
{{\mathscr Z}(n,x) \over n} \equiv
\!\!\! \sum_{0 \le \ell \le n/2} \! (-1)^\ell {n \!-\! \ell \choose \ell} 
(2x \!-\! 1)^{n-2\ell} \, {{\mathcal Z}(n \!-\! \ell,v) \over n-\ell} \quad 
\mbox{for } n=1,2,\ldots , \nonumber
\end{equation}
which clearly generalize eqs.(\ref{Z10}) (for $n$ odd) 
and (\ref{Z1E}) (for $n$ even) away from $x=\hf$.
On the other hand, eqs.(\ref{LID}) invert
into finite linear relations of the same form (if $x \ne \hf$):
$ (\tau^2+v)^{-m} \equiv \sum\limits_{n=1}^m V_{m,n}(x) 
[(x-\rho)^{-n} + (x-1+\rho)^{-n}] $. At this point we can identify $V_{m,n}(x)$:
it has to be the coefficient of $(x-\rho)^{-n}$ in the Laurent series of 
$ (\tau^2+v)^{-m} = (x-\rho)^{-m} [2x-1-(x-\rho)]^{-m} $ 
in powers of $(x-\rho)$, i.e., $V_{m,n}(x)={2m-n-1 \choose m-1} (2x-1)^{n-2m}$.
Then, upon summing the above expression for $(\tau^2+v)^{-m}$ 
over half the zeros, we finally get
\begin{eqnarray}
\label{ZIL}
{\mathcal Z}(m,v) \si \equiv \sf \sum_{\ell=0}^{m-1} 
{m \!+\! \ell \!-\! 1 \choose m \!-\! 1} (2x \!-\! 1)^{-m-\ell} 
{\mathscr Z}(m \!-\! \ell,x) \qquad \qquad \\ 
&& \qquad \qquad \qquad \qquad 
\mbox{for } v \equiv (x \!-\! \hf)^2 \ne 0 \mbox{ and } m=1,2,\ldots , \nonumber
\end{eqnarray}
whereas ${\mathcal Z}(m,0) \equiv \hf (-1)^m {\mathscr Z}(2m,\hf)$, 
by eq.(\ref{Z1E}). 
Remark: the pair of mutually inverse relations (\ref{LIZ}) and (\ref{ZIL}) 
are clearly similar to the identities (\ref{ZSK}) and have the same origin.
They extend to all primary functions $L$ and all $v$-values 
previous results written only for the Riemann case and $v={1 \over 4}$ 
\cite{Ma,Vz}.
 
\begin{table}
\centering
\begin{tabular} {cc}
\hline \\[-12pt]
$\sigma$ & 
$ {\mathcal Z}(\sigma,v) = \sum\limits_{k=1}^\infty ({\tau_k}^2+v)^{-\sigma} $ 
\\[10pt]
\hline \\[-12pt]
$-m \le 0$ & $ \hf \Bigl[ -\sum\limits_{j=0}^m \, 
{\textstyle m_{\vphantom{o}} \choose \textstyle j^{\vphantom{0}}}
(-1)^j {\mathbf Z}(-2j,\hf) v^{m-j} + q {(v-{1 \over 4})}^m \Bigr] $ \\
$0$ & $ \hf(\hf \tilde a_1 + \tilde a_0 +q)$ \\[4pt]
\hline \\[-12pt]
${\mbox{\sl $\sigma$-derivative} \atop \mbox{\sl at 0}}$ & 
$ {\mathcal Z}'(0,v) = \hf \, b_1 + b_0 - \log \Xi(\hf \pm v^{1/2}) $ \\[6pt]
$+m \ge 1$ & $\!\!\!\!\!\!\!\!\! \left\{ \matrix{
\!\! - \! \sum\limits_{\ell=0}^{m-1}
{\textstyle m \!+\! \ell \!-\! 1_{\vphantom{0}} \choose 
\textstyle m \!-\! 1^{\vphantom{0}}} \,
{\textstyle (\mp 2 v^{1/2})^{-m-\ell} \over \textstyle (m \!-\! \ell \!-\! 1)!} 
(\log \Xi)^{(m-\ell)}(x)|_{x={1 \over 2} \pm v^{1/2} } &
(v \ne 0) \cr
{\textstyle (-1)^{m+1} \over \textstyle 2(2m-1)!} (\log \Xi)^{(2m)}(\hf) \hfill 
& (v=0) \cr} \right. $ \\[24pt]
\hline\\[-12pt]
\end{tabular}
\caption{Special values of ${\mathcal Z}(\sigma,v)$
(upper half: algebraic, lower half: transcendental \cite{Vz})
for a general primary zeta function $L(x)$ with a pole of order $q$ (at $x=1$).
Notations: see eqs.(\ref{FE}) for $\Xi(x)$,
(\ref{SDef}) and (\ref{TI}) for ${\mathbf Z}(-n,x)$,
(\ref{ST}) for $\tilde a_j,\ b_j$; $m$ is integer.
}
\end{table}

\medskip

The resulting values of ${\mathcal Z} (\sigma,v)$ for general $v$
are listed in Table~2 (lower half). 
(This Table improves upon Table~1 of \cite{VZ} in two independent ways: 
it is valid for zeros of a general primary function $L(x)$, not just $\zeta(x)$,
and it specifies ${\mathcal Z}(+m,v)$ more explicitly.)

In the particular case $v={1 \over 4}$,
the transcendental special values of ${\mathcal Z}(\sigma,{1 \over 4})$ 
involve those of ${\mathscr Z}(s,1)$ by eq.(\ref{ZIL}), hence they will likewise
end up expressed in terms of the generalized Stieltjes cumulants (\ref{CUM}), 
cf. Tables~3, 4, 6 below.

\section{The third family $\{ {\mathfrak Z}(\sigma,y) \}$}

As with the preceding case, a starting point is the knowledge of 
one particular member of the family, now through the obvious identity 
${\mathfrak Z}(\sigma,0) \equiv {\mathcal Z}(\sigma,0)$.
All results of Sec.~3.1 then cover this case as well.

The generic ${\mathfrak Z}(\sigma,y)$ (with $y \ne 0$) is built on 
a {\sl desymmetrized\/} set of zeros, say $(\hf + {\rm i} \tau_k)$ only,
hence it will be harder to describe explicitly than the other two families.
Still, its polar structure can be drawn directly from an expansion 
(in $\{ |y| <\min_k \{ |\tau_k| \} \}$) 
similar to eq.(\ref{SHZ0}) for ${\mathcal Z} (\sigma,v)$ (see also \cite{HKW}):
\begin{equation}
\label{SHZ1}
{\mathfrak Z}(\sigma,y) = \sum_{k=0}^\infty 
{\tau_k}^{-2\sigma} \Bigl( 1 + {y \over \tau_k} \Bigr)^{-2\sigma} = 
\sum_{\ell=0}^\infty 
{\Gamma(1-2\sigma) \over \ell ! \Gamma(1 \!-\! 2\sigma \!-\! \ell)}
{\mathcal Z}(\sigma + \hf \ell) \, y^\ell .
\end{equation}
This formula generates a pole for ${\mathfrak Z}(\sigma,y)$
now at every {\sl integer or half-integer\/} $\hf(1-n)$, $n \in\mathbb N$, 
according to:
\begin{equation}
\label{HP1}
{\mathfrak Z}(\hf(1-n)+\varepsilon,y) = \sum_{\ell=0}^n 
{\Gamma(n-2\varepsilon) \over \ell ! \Gamma(n \!-\! \ell \!-\! 2\varepsilon)}
{\mathcal Z}(\hf(1-n+\ell)+\varepsilon) \, y^\ell \ + \Biggl\{ \matrix
{ O(1) \hfill &\mbox{for }n=0 \hfill \cr 
O(\varepsilon) &\mbox{for }n= 1,2,\ldots \cr}
\Biggr. 
\end{equation}
Concrete differences with eq.(\ref{HP2}) arise from the factor
$\Gamma(n-2\varepsilon)/\Gamma(n-\ell-2\varepsilon)$ vanishing 
whenever $\ell \ge n >0$. Only the polar part at $\sigma=\hf$ remains 
the same as for ${\mathcal Z}(\sigma,v)$
(of order $r=2$ and independent of $y$, given by eq.(\ref{Z2P})); 
all other poles $\hf(1-n)$ of ${\mathfrak Z}(\sigma,y)$ are now {\sl simple\/}, 
of residues
\begin{equation}
\label{RES1}
r_n(y) = -{\tilde a_1 \over 2 \pi n} \, y^n +
\sum_{0< 2m \le n} { n-1 \choose 2m-1 } {\mathcal R}_m \, y^{n-2m},
\quad n=1,2,\ldots 
\end{equation}
(in terms of the residues ${\mathcal R}_m$ given by eq.(\ref{Z2Q})).
Moreover, at $\sigma=0$ (only), the $\varepsilon$-expansion of eq.(\ref{HP1}) 
captures the {\sl finite part\/} too (cf. eqs.(\ref{ST}),(\ref{Z0V})):
\begin{equation}
\label{RES0}
r_1(y) = {\rm Res}_{\sigma=0} {\mathfrak Z}(\sigma,y) 
= -{\tilde a_1 \over 2 \pi} \, y \, ;
\quad {\rm FP}_{\sigma=0} \,{\mathfrak Z}(\sigma,y) = 
{\textstyle{1 \over 4}} \tilde a_1 + \hf(\tilde a_0 +q) - { b_1 \over \pi } \, y .
\end{equation}

On the other hand, while previously we could express infinite products 
over {\sl all\/} the zeros such as $\mathscr D$, resp. $\mathcal D$,
in terms of simpler functions like $\mathbf D$ and $L$ 
(through eqs.(\ref{zf1}), resp. (\ref{zf2})), 
now we lack that ability for a similar infinite product 
but restricted to {\sl half\/} the zeros.
We thus have no simple formulae for transcendental special values of 
${\mathfrak Z}(\sigma,y)$ (i.e., in the half-plane $\{ {\rm Re\,} \sigma >0 \}$).
Only a sequence of binary relations results by specializing 
the identity (\ref{Z2X}) to $s \in \mathbb N^\ast$, 
\begin{equation}
\bigl[ {\rm i}^m {\mathfrak Z} \bigl( \hf m, {\rm i} t \bigr) +
{\rm i}^{-m} {\mathfrak Z} \bigl( \hf m, -{\rm i} t \bigr) \bigr] 
\equiv {\mathscr Z} \bigl( m,\hf+t \bigr) , \qquad m=1,2,\ldots,
\end{equation}
constituting an obvious result except for $m=1$: 
then, {\sl finite parts\/} are to be taken {\sl on the left-hand side only\/}.

\section{Concrete examples}

We finally illustrate the preceding results upon the two classes 
of primary zeta functions highlighted in Sec.~1.3. 
As a rule, it suffices to specialize the general formulae as indicated below
(and recalling that the order $\mu_0=1$ throughout).
To strengthen the practical side of this work,
we will further display the final concrete formulae reached 
when the shift parameters are themselves fixed at specially interesting values.
We will mainly show results for ${\mathscr Z}(s,1) \ [\equiv {\mathscr Z}(s)]$ 
and ${\mathscr Z}(s,\hf)$, as Tables~4--7
(generalizing Table~3 in \cite{VZ}, which was specific to $L(x)=\zeta(x)$).
The analogous special cases for the second family, ${\mathcal Z}(\sigma,0)$ 
$[\equiv {\mathfrak Z}(\sigma,0) \equiv {\mathcal Z}(\sigma) ]$
and ${\mathcal Z}(\sigma,{1 \over 4}) $, 
are easily recovered from the preceding ones by applying the results of Sec.~3,
specially as subsumed in Table~3;
this stage entails an extension of Table~1 in \cite{Vz} to general $L(x)$.
(We do not further illustrate the third family $\{ {\mathfrak Z}(\sigma,y) \}$,
for its lack of explicit special values.)

\begin{table}
\centering
\begin{tabular} {ccc}
\hline \\[-12pt]
$\sigma$ & ${\mathcal Z}(\sigma,0) = 
\sum\limits_{k=1}^\infty {\tau_k}^{-2\sigma}$ &
$ {\mathcal Z}(\sigma,{1 \over 4}) = 
\sum\limits_{k=1}^\infty {({\tau_k}^2+{1 \over 4})}^{-\sigma} $ 
\\[10pt]
\hline \\[-12pt]
$-m \le 0$ & $ \hf (-1)^m {\mathscr Z}(-2m,\hf) $ & $ \hf \sum\limits_{j=0}^m 
\, {\textstyle m_{\vphantom{o}} \choose \textstyle j^{\vphantom{0}}}
(-1)^j 2^{-2(m-j)} {\mathscr Z}(-2j,\hf) $ \\
$0$ & $ \hf {\mathscr Z}(0,\hf) $ & $ \hf {\mathscr Z}(0,\hf) $ \\[4pt]
\hline \\[-12pt]
${\mbox{\sl derivative} \atop \mbox{\sl at 0}}$ & 
$ {\mathcal Z}'(0,0) = {\mathscr Z}'(0,\hf) $ &
$ {\mathcal Z}'(0,{1 \over 4}) = -\hf \, b_1 + {\mathscr Z}'(0,1) $ \\[2pt]
$+m \ge 1$ & $ \hf (-1)^m {\mathscr Z}(2m,\hf) $ & 
$ \sum\limits_{\ell=0}^{m-1}
{\textstyle m \!+\! \ell \!-\! 1_{\vphantom{0}} \choose 
\textstyle m \!-\! 1^{\vphantom{0}}} {\mathscr Z}(m-\ell,1) $ \\[10pt]
\hline\\[-12pt]
\end{tabular}
\caption{General formulae for the special values of 
${\mathcal Z}(\sigma,0) \equiv {\mathcal Z}(\sigma)$ 
and ${\mathcal Z}(\sigma,{1 \over 4})$, 
in terms of those for ${\mathscr Z}(s,\hf)$ and ${\mathscr Z}(s,1)$ 
(such as provided in the subsequent Tables).
Notations: see eq.(\ref{ST}) for $b_1$; 
$m$ is integer.
}
\end{table}

At those special parameter cases, 
the relevant values of $\mathbf Z(s,x)$ will also become more explicit, using
\begin{equation}
\zeta(s,1) \equiv \zeta(s); \quad \zeta(s,\hf) \equiv (2^s-1)\zeta(s); \quad
2^{-2s} \zeta(s,\hf \mp {\textstyle{1 \over 4}}) 
\equiv \hf[(1-2^{-s}) \zeta(s) \pm \beta(s)] .
\end{equation}
(cf. eq.(\ref{BET})).
So, the Dirichlet series $\zeta(s)$ and $\beta(s)$ ($\equiv L_{\chi_4}(s)$;
$\chi_4$ is the real primitive character for the modulus 4)
will both occur ubiquitously in these particular special values,
for whatever choice of primary Dirichlet series $L(x)$.

Again, the Riemann case $L(x) = \zeta(x)$, having $q=1$,
is more conveniently treated here as a special Dedekind zeta function only
(not an $L$-function).

\subsection{$L$-functions for real primitive Dirichlet characters}

According to eq.(\ref{LV}), any $L$-function $L_\chi(x)$ 
for such a Dirichlet character $\chi \pmod d$
(with $d>1$, to exclude $\zeta(x)$) is handled by the choice
\begin{equation}
\label{VL}
q \equiv 0; \qquad 
{\mathbf G}(x) \equiv (\pi /d)^{x/2} / \Gamma \bigl( \hf(x+a) \bigr) ,
\quad a= \biggl\{ \matrix{ 0 & \mbox{ for $\chi$ even} \cr 
1 & \mbox{ for $\chi$ odd.} \hfill \cr} \Biggr.
\end{equation}

In turn, the other useful quantities specialize as follows:

\noindent - the leading coefficients in the Stirling formula (\ref{ST}):
\begin{equation}
\matrix{\tilde a_1=\hf \, , \hfill && \tilde a_0=\hf(a-1) , \hfill \cr \cr
b_1=-\hf \log(2\pi/d) , \hfill && b_0=\hf \log(2^{2-a} \pi) ; \hfill \cr} 
\end{equation}
- the shadow zeta function (eq.(\ref{SDef})):
\begin{equation}
{\mathbf Z}(s,x) = 2^{-s} \zeta \bigl( s,\hf (x+a) \bigr) ;
\end{equation}
- the lowest generalized Stieltjes cumulants (eq.(\ref{CUM})): 
$g_0^{\rm c} \{ L_\chi \} \equiv -\log L_\chi(1)$ can always be specified,
as well as $g_1^{\rm c} \{ L_\chi \} \equiv [{L_\chi}'/L_\chi](1)$ when $a=1$. 
First, the general formula 
\begin{equation}
\label{LHZ}
L_\chi(x) \equiv d^{-x} \sum\limits_{n=1}^d \chi(n) \, \zeta(x,n/d),
\end{equation}
together with the special values (\ref{HN}--\ref{HD}) 
of the Hurwitz zeta function
(also using $\chi(d)=0$ and $\sum\limits_{n=1}^d \chi(n)=0$ throughout),
yield these special values for $L_\chi(x)$:
\begin{eqnarray}
\label{L0}
L_\chi(0) \si = \sf -{1 \over d} \sum_{n=1}^{d-1} \chi(n) \, n 
\qquad \qquad \qquad \qquad \qquad \qquad \mbox{(algebraic)} \\
\label{L1}
L_\chi(1) \si = \sf 
-{1 \over d} \sum_{n=1}^{d-1} \chi(n) {\Gamma'(n/d) \over \Gamma(n/d)} 
\qquad \qquad \qquad \qquad \mbox{(transcendental)} \\
\label{LP0}
{L_\chi}'(0) \si = \sf 
-L_\chi(0) \log d + \sum_{n=1}^{d-1} \chi(n) \log \Gamma(n/d)
\qquad \mbox{(transcendental)} .
\end{eqnarray}
Then, the functional equation (\ref{FE}) also implies 
(we now suppress the $\chi$ labels)
\begin{eqnarray}
\!\!\!\!\!\!\!\! \mbox{if }a=1: \ \ L(1) \si = \sf \pi d^{-1/2} L(0), 
\ \ [L'/L](1) = \gamma +\log(2\pi/d) - [L'/L](0) ;\\
\!\!\!\!\!\!\!\! \mbox{if }a=0: \ \ L(1) \si = \sf 2 d^{-1/2} L'(0) \quad
\mbox{[$L'(1)$ involves the unknown $L''(0)$,...]}
\end{eqnarray}
(in the $a=0$ case, moreover, $L(0) \equiv 0$, and eq.(\ref{LP0}) for $L'(0)$ 
simplifies further by the reflection formula for $\Gamma$, 
with $\chi$ being even). The final outcome is:

\begin{table}
\centering
\begin{tabular} {ccc}
\hline \\[-7pt]
$s$ & &
${\mathscr Z}(s) \defi {\mathscr Z}(s,1) \equiv 
\sum\limits_\rho \rho^{-s} \qquad [x=1]$ \\[10pt]
\hline \\[-10pt]
$-n \!<\! 0$ & & $ \Bigl[ (a \!-\! 1)(2^n \!-\! 1) + a \, 2^n \Bigr] 
{\textstyle B_{n+1} \over \textstyle n \!+\! 1} $ \\[4pt]
$0$ & & $\hf a$ \\[4pt]
\hline \\[-12pt]
{\sl derivative at 0} & &
$ {\mathscr Z}'(0) = 
\hf \bigl[ (1 \!-\! a) \log 2 + a \log \pi \bigr] + g_0^{\rm c} \{ L_\chi \} $ 
\\[4pt]
{\sl finite part at +1} & &
${\rm FP}_{s=1} {\mathscr Z}(s) = 
(a \!-\! \hf)\log 2 - \hf \gamma + g_1^{\rm c} \{ L_\chi \}$
\\[6pt]
$+1$ & & 
$ (a \!-\! 1)\log 2 - \hf \log (\pi/d) - \hf \gamma + g_1^{\rm c} \{ L_\chi \} $ 
\\[4pt]
$+n \!>\! 1$ & & $ [(a \!-\! 1) (1 \!-\! 2^{-n}) - a \, 2^{-n}] \, \zeta(n)
+ {\textstyle g_n^{\rm c} \{ L_\chi \} \over \textstyle (n \!-\! 1)!} $ \\[8pt]
\hline\\[-12pt]
\end{tabular}
\caption{Special values of the zeta function ${\mathscr Z}(s,1)$ 
over the nontrivial zeros of an $L$-function for 
a real primitive Dirichlet character $\chi$ 
of modulus $d>1$ and parity $a=0$ or 1 (see eqs.(\ref{KH}--\ref{A01})). 
For the $g_n^{\rm c}$, see eqs.(\ref{CUM}), (\ref{GO}--\ref{C34}).
Notations: see eqs.(\ref{NOT}); $n$ is integer.
In last line, $\zeta(n) \equiv (2\pi)^n |B_n|/(2 \, n!)$ when $n$ is even.
}
\vskip 1cm

\begin{tabular} {cc}
\hline \\[-7pt]
$s$ &
${\mathscr Z}(s,\hf) \equiv \sum\limits_\rho {(\rho-\hf)}^{-s} \qquad [x=\hf]$ 
\\[10pt]
\hline \\[-12pt]
even $-n \! \le \! 0$ & $ 2^{-n-1}(a \!-\! \hf) E_n $ \\ [4pt]
~odd $-n \!<\! 0$ &
$ -\hf (1 \!-\! 2^{-n}){\textstyle B_{n+1} \over \textstyle n \!+\! 1} $ 
\\[4pt]
$0$ & $\hf (a \!-\! \hf)$ \\[4pt]
\hline \\[-12pt]
{\sl derivative at 0} &
$ {\mathscr Z}'(0,\hf) = ({3 \over 4} \!-\! a) \log 2 
+ (a \!-\! \hf) \log \bigl[ \Gamma({1 \over 4})^2/\pi \bigr] - \log L_\chi (\hf) $ 
\\[4pt]
{\sl finite part at +1} &
$ {\rm FP}_{s=1} {\mathscr Z}(s,\hf) = \hf \log (2 \pi /d) $ \\[6pt]
~odd $+n \! \ge \! 1$ & 0 \\[-4pt]
even $+n \!>\! 1$ & 
$-\hf \bigl[ (2^n \!-\! 1) \,\zeta (n) + (1 \!-\! 2a) \, 2^{n} \beta (n) \bigr]
- {\textstyle (\log L_\chi)^{(n)} (\hf) \over \textstyle (n \!-\! 1)!} $ 
\\[8pt]
\hline\\[-12pt]
\end{tabular}
\caption{Same as above, but for the zeta function ${\mathscr Z}(s,\hf)$.
In last line, $n$ being even, $\zeta(n) \equiv (2\pi)^n |B_n|/(2 \, n!)$
while $\beta(n)$ (see eq.(\ref{BET})) remains elusive.
}
\end{table}

\noindent $\bullet$ when $a=1$, an algebraic explicit formula for $L(1)$ 
plus a transcendental one for $[L'/L](1)$ (in terms of Gamma values), 
amounting to
\begin{equation}
\label{GO}
\!\! g_0^{\rm c} \{ L_\chi \} = 
-\log \Bigl[ -{\pi \over d^{3/2}} \sum_{n=1}^{d-1} \chi(n)\, n\Bigl],
\ \ g_1^{\rm c} \{ L_\chi \} = \gamma +\log (2\pi) + { \sum\limits_{n=1}^{d-1} 
\chi(n) \log \Gamma \Bigl( \displaystyle{n \over d} \Bigr) 
\over \sum\limits_{n=1}^{d-1} \chi(n) \, \displaystyle{n \over d} }\, ;
\end{equation}
$\bullet$ when $a=0$, just a formula for $L(1)$,
transcendental but more elementary than (\ref{L1}), and amounting to
\begin{equation}
\label{GE}
\!\!\!\!\!\! g_0^{\rm c} \{ L_\chi \} = -\log \Bigl[ -{1 \over d^{1/2}} 
\sum_{n=1}^{d-1} \chi(n) \log \sin {\pi n \over d} \Bigr] 
\quad \qquad
[g_1^{\rm c} \{ L_\chi \} \mbox{ not specified]} .
\end{equation}
E.g., for each of $d=3$ and 4 (the lowest possible values of $d$),
the real primitive character $\chi_d$ is unique
($\chi_d (\pm 1 \bmod d) = \pm 1$, else $\chi_d(n)=0$;
in particular, $L_{\chi_4}(x) \equiv \beta(x)$ as in eq.(\ref{BET}));
they are both odd, giving
\begin{equation}
\!\label{C34}
\matrix{ 
d=3: & g_0^{\rm c} \{ L_{\chi_3} \} \!\si = \sf\! -\log (\pi/3^{3/2}), 
& g_1^{\rm c} \{ L_{\chi_3} \} \!\si = \sf\! 
\log [(2\pi)^4/3^{3/2}] + \gamma - 6 \log \Gamma({1 \over 3}) ; \cr
d=4: & \hfill g_0^{\rm c} \{ L_{\chi_4} \} \!\si = \sf\! -\log (\pi/4), \hfill 
& \hfill g_1^{\rm c} \{ L_{\chi_4} \} \!\si = \sf\! 
\log(4 \pi^3) + \gamma -4 \log \Gamma({1 \over 4}) . \hfill \cr}
\end{equation}
In general, whether $a=0$ or 1, we cannot specify the $g_n^{\rm c}$ any further
than stated; we might just relate them, through eq.(\ref{LHZ}), 
to special cases of the still more general Laurent coefficients 
$\gamma_n(w)$ of $\zeta(x,w)$ around $x=1$ \cite{BE,Kr}.
\medskip

The finally resulting special values of ${\mathscr Z}(s,1)$ and 
${\mathscr Z}(s,\hf)$, over zeros of a general real primitive 
Dirichlet character, are presented in Tables~4 and 5 respectively.

\subsection{Dedekind zeta functions}

Referring back to eq.(\ref{DV}), the Dedekind zeta function $\zeta_K(x)$
of any algebraic number field $K$ is handled by the choice
\begin{equation}
\label{VD}
q \equiv 1, \qquad {\mathbf G}(x) \equiv 
{ \bigl( 4^{r_2} \pi^{n_K} / |d_K| \bigr)^{x/2} \over 
x \Gamma (x/2)^{r_1} \Gamma (x)^{r_2} } \, .
\end{equation}
Moreover, $r \defi {\rm Res}_{x=1} \zeta_K(x)$ (the sole residue of $\zeta_K(x)$) 
is strictly positive and {\sl computable\/}
(in terms of many field invariants) (\cite{He} thms~121, 124; \cite{St}).

In turn, the other useful quantities specialize as follows:

\begin{table}
\centering
\begin{tabular} {ccc}
\hline \\[-7pt]
$s$ & &
${\mathscr Z}(s) \defi {\mathscr Z}(s,1) \equiv
\sum\limits_\rho \rho^{-s} \qquad [x=1]$ \\[10pt]
\hline \\[-10pt]
$-n \!<\! 0$ & & $ \bigl[-r_1(2^n \!-\! 1) + r_2 \bigr] 
\, {\textstyle B_{n+1} \over \textstyle n \!+\! 1} + 1 $ \\[4pt]
$0$ & & $\hf r_2 + 2$ \\[4pt]
\hline \\[-12pt]
{\sl derivative at 0} & &
$ {\mathscr Z}'(0) = \hf \bigl[ (r_1 \!+\! r_2) \log 2 + r_2 \log \pi \bigr]
+ g_0^{\rm c} \{ \zeta_K \} $ \\[4pt]
{\sl finite part at +1} & &
${\rm FP}_{s=1} {\mathscr Z}(s) = - \hf r_1 \log 2 
+ 1 - \hf n_K \gamma + g_1^{\rm c} \{ \zeta_K \} $ \\[6pt]
$+1$ & & 
$ \hf \log |d_K| \!-\! (r_1 \!+\! r_2) \log 2 \!-\! \hf n_K \log \pi 
\!+\! 1 \!-\! \hf n_K \gamma + g_1^{\rm c} \{ \zeta_K \} $ \\[4pt]
$+n \!>\! 1$ & & $ - [r_1(1 \!-\! 2^{-n}) \!+\! r_2] \, \zeta(n) + 1
+ {\textstyle g_n^{\rm c} \{ \zeta_K \} \over \textstyle (n \!-\! 1)!} $ \\[8pt]
\hline\\[-12pt]
\end{tabular}
\caption{Special values of the zeta function ${\mathscr Z}(s,1)$ 
over the nontrivial zeros of a Dedekind zeta function 
for an algebraic number field $K$ (see eqs.(\ref{DZ}--\ref{DV})).
For the $g_n^{\rm c}$, see eqs.(\ref{CUM}), (\ref{CZ}--\ref{CQ}).
Notations: see eqs.(\ref{NOT}); $n$ is integer.
In last line, $\zeta(n) \equiv (2\pi)^n |B_n|/(2 \, n!)$ when $n$ is even.
}
\vskip 1cm

\begin{tabular} {cc}
\hline \\[-7pt]
$s$ &
${\mathscr Z}(s,\hf) \equiv \sum\limits_\rho {(\rho-\hf)}^{-s} \qquad [x=\hf]$ 
\\[10pt]
\hline \\[-12pt]
even $-n \! \le \! 0$ & $ 2^{-n+1}(1-{1 \over 8} r_1 E_n) $ \\ [4pt]
~odd $-n \!<\! 0$ &
$ -\hf n_K (1 \!-\! 2^{-n}) \, {\textstyle B_{n+1} \over \textstyle n \!+\! 1} $
\\[4pt]
$0$ & $2 - {1 \over 4} r_1$ \\[4pt]
\hline \\[-12pt]
{\sl derivative at 0} &
$\!\! {\mathscr Z}'(0,\hf) = (2 \!+\! {3 \over 4} r_1 \!+\! \hf r_2) \log 2 
\!-\! \hf r_1 \log \bigl[ \Gamma({1 \over 4})^2/\pi \bigr] 
\!-\! \log |\zeta_K| (\hf) $ 
\\[4pt]
{\sl finite part at +1} &
$ {\rm FP}_{s=1} {\mathscr Z}(s,\hf) = 
\hf \bigl[ n_K \log (2 \pi) - \log |d_K| \bigl] $ 
\\[6pt]
~odd $+n \! \ge \! 1$ & 0 \\
even $+n \!>\! 1$ & 
$-\hf n_K (2^n \!-\! 1) \,\zeta (n) \!-\! \hf r_1 \, 2^{n} \beta (n) 
\!+\! 2^{n+1} \!-\! 
{\textstyle (\log |\zeta_K|)^{(n)} (\hf) \over \textstyle (n \!-\! 1)!} $
\\[8pt]
\hline\\[-12pt]
\end{tabular}
\caption{Same as above, but for the zeta function ${\mathscr Z}(s,\hf)$.
In last line, $n$ being even, $\zeta(n) \equiv (2\pi)^n |B_n|/(2 \, n!)$
while $\beta(n)$ (see eq.(\ref{BET})) remains elusive.
}
\end{table}

\noindent - the leading coefficients in the Stirling formula (\ref{ST}):
\begin{equation}
\matrix{\tilde a_1=\hf n_K \, , \hfill && \tilde a_0=1-\hf(r_1+r_2) , 
\hfill \cr\cr
b_1=-\hf \bigl[ \log \bigl( (2\pi)^{n_K}/|d_K| \bigr) \bigr] , \hfill && 
b_0=(r_1 \!+\! \hf r_2) \log 2 + \hf (r_1 \!+\! r_2) \log \pi ; \hfill \cr} 
\end{equation}
- the shadow zeta function (eq.(\ref{SDef})): 
counting all zeros of $\mathbf G(x)$ with their multiplicities, it reads as
\begin{equation}
{\mathbf Z}(s,x) = r_1 2^{-s} \zeta ( s,\hf x ) + r_2 \zeta(s,x) -x^{-s} ;
\end{equation}
- the lowest generalized Stieltjes cumulants (eq.(\ref{CUM})): 
$g_0^{\rm c} \{ \zeta_K \} \equiv -\log r$ is fully explicit and, 
sometimes at least, 
$g_1^{\rm c} \{ \zeta_K \} \equiv r^{-1} \, {\rm FP}_{x=1} \zeta_K(x)$ 
can also be described. 
First, if $K= \mathbb Q$: then $\zeta_K(x) \equiv \zeta(x)$ with $r=1$, hence
\begin{equation}
\label{CZ}
g_0^{\rm c} \{ \zeta \} = 0, \qquad \qquad g_1^{\rm c} \{ \zeta \} = \gamma ;
\end{equation}
in turn, a general $g_1^{\rm c}$ will be some extension of 
Euler's constant $\gamma$.
As next example, if $K$ is a {\sl quadratic\/} number field, 
then $\zeta_K(x) \equiv \zeta(x) \, L_\chi(x)$
where $\chi$ is the real primitive character of modulus $|d_K|$ 
given by the Kronecker symbol for the discriminant $d_K$ (\cite{He} Sec.~49).
Now in general, 
the zeta functions over the zeros (and their linear invariants) obviously
{\sl add up\/} when their primary functions $L$ are {\sl multiplied\/}.
So, for a quadratic number field,
\begin{equation}
\label{CQ}
g_0^{\rm c} \{ \zeta_K \} = g_0^{\rm c} \{ L_\chi \}, \qquad 
g_1^{\rm c} \{ \zeta_K \} = \gamma + g_1^{\rm c} \{ L_\chi \}
\end{equation}
referring to the same cumulants for $L_\chi(x)$
that were precisely described under the previous example 
by eq.(\ref{GO}) for $\chi$ odd, or (\ref{GE}) for $\chi$ even.
Two basic examples (both with $r_1=0$, $r_2=1$) are: $K= \mathbb Q({\rm i})$ 
(for which $d_K=-4$, $\chi = \chi_4$,
$L_\chi(x) \equiv \beta(x)$ as in eq.(\ref{BET})),
and $K= \mathbb Q(\sqrt{-3})$ (for which $d_K=-3$, $\chi = \chi_3$),
hence their specific cumulants $g_0^{\rm c}$, $g_1^{\rm c}$ 
follow from eqs.(\ref{C34}), (\ref{CQ}).
\medskip

The finally resulting special values of ${\mathscr Z}(s,1)$ and 
${\mathscr Z}(s,\hf)$, over zeros of a general Dedekind zeta function, 
are presented in Tables~6 and 7 respectively.

\end{document}